\documentclass[12pt]{amsart}
\usepackage{amsmath}
\usepackage{amsfonts,amssymb,amsthm, txfonts, pxfonts, boxedminipage}

\DeclareMathOperator{\acosh}{arccosh}

\newtheorem{theorem}{Theorem}
\newtheorem{conjecture}{Conjecture}
\newtheorem{lemma}{Lemma}
\newtheorem{corollary}{Corollary}
\theoremstyle{remark}
\newtheorem{remark}[theorem]{Remark}

\newtheorem{observation}[theorem]{Observation}
\begin{document}

\title[Volumes and degeneration]{Volumes of degenerating polyhedra -- on
  a conjecture of J.~W.~Milnor}
\author{Igor Rivin}
\begin{abstract}
In his paper \cite{milnorcol} J. Milnor conjectured that the volume $V_n$ of
$n$-dimensional hyperbolic and spherical simplices, as a function of the
dihedral angles, extends continuously to the closure
$\overline{\mathbb{A}}$ of the space $\mathbb{A}$ of
allowable angles (``The continuity conjecture'') , and furthermore, $V_n(a\in \partial \mathbb{A}) = 0$
if and only if $a$ lies in the closure of the space of angles of
Euclidean tetrahedra (``the Vanishing Conjecture'').
A proof of the Continuity Conjecture  was given by F. Luo
(\cite{luovol}-- Luo's argument uses Kneser's formula \cite{kneser36}
  together with some delicate geometric estimates). 
In this paper we give a simple proof of both parts of Milnor's conjecture, prove much sharper regularity results, and then extend the method to apply to all 
convex polytopes. We also give a precise description of 
the boundary of the space of angles of convex polyhedra in
$\mathbb{H}^3,$ and sharp estimates on the diameter of a polyhedron in
terms of the length of the shortest polar geodesic.
\end{abstract} 
\thanks{The author would like to thank D. Laksov for his continued encouragement.}

\date{\today}

\maketitle
\section{Introduction}
Consider the set of simplices in $\mathbb{H}^n$ or $\mathbb{S}^n.$ It is
well-known that this set is parametrized by the (ordered) collection of
dihedral angles, and we may call the set of assignements of dihedral angles of
geometric simplices in $\mathbb{H}^n$ as a subset $\Omega_{\mathbb{H}^n}
\subset \mathbb{R}^{n(n+1)/2 }$ and similarly, the set of dihedral angle
assignements of of geometric simplices in $\mathbb{S}^n$ as a subset
$\Omega_{\mathbb{S}^n} \subset \mathbb{R}^{n(n+1)/2}.$  These sets are open,
since they are defined by collections of strict inequalities (which are
polynomial in the cosines of the dihedral angles).  One may then view the
volume $V$ of a simplex as a function $V$ on $\Omega_{X^n}.$ J.~Milnor
(\cite{milnorcol})
conjectured:

\begin{conjecture}
\label{milnorconj}
The volume function $V$ admits a continuous extension to
$\overline{\Omega}.$ Furthermore, the points on $\partial{\Omega}$
where $V$ vanishes are precisely those which also lie in the closure
of the set of angle assignments of Euclidean simplices.
\end{conjecture}

Some comments are in order regarding Conjecture \ref{milnorconj}.
Firstly, it falls into two parts: the ``Continuity Conjecture'' and
the ``Vanishing Conjecture.'' The Vanishing Conjecture cannot be
stated without knowing that the answer to the Continuity Conjecture is
affirmative.

The Continuity Conjecture was first shown by F.~Luo (in
\cite{luovol}), and then a sharper version was shown by me in a
a predecessor (\cite{rivvol}) of the current paper. 

Milnor does not attribute the conjecture to himself, and his paper
(which was written in the late seventies or early eighties) seems to
imply that the conjecture precedes the paper. 

The contents of this paper are as follows. First, we give a simple
argument to show a sharp version of Milnor's Continuity Conjecture for
all hyperbolic polytopes of dimension greater than $3,$ and also all
spherical polytopes. It should be noted that since in many cases it is
not known whether hyperbolic or spherical polytopes are determined by
their dihedral angles and how to characterize the possible assignments
of dihedral angles\footnote{Simplices are a notable exception, and an
excellent exposition is given in Milnor's paper \cite{milnorcol}}, it
makes more sense to use \emph{polar metrics} introduced in
\cite{thes0, thes}. The argument shows that the extension is, in
fact, Lipschitz. 

Next we give an argument to show the Continuity Conjecture for
three-dimensional hyperbolic tetrahedra, which is conceptually related
to the higher-dimensional argument (\emph{via} the Schl\"afli
differential formula), but is a little more delicate. The argument
requires a version of Sobolev's Embedding Theorem, but as a
consequence, a sharp regularity result is obtained (this time the
extension is shown to be in the class $C^{0, 1}.$)

We then go on to arbitrary convex polyhedra in $\mathbb{H}^3$ (and
polytopes in $\mathbb{H}^n$) and
prove the same sharp version of the Continuity Conjecture for
those. These results use (at least philosophically) the results of
\cite{thes0,thes}.It should be noted that the estimates proved in
this section work just as well for higher-dimensional convex polytopes
(although they are not necessary for the regularity result). The
results here are of independent interest, and can be summarized as
follows:
\begin{theorem}
Let $P$ be a polyhedron with $N$ vertices in $\mathbb{H}^3$ of diameter $\rho
\gg 1.$ Let $M^*$ be the polar metric of $P$ (as in \cite{thes,thes0}). The
$M^*$ lies within $c_1(N) \exp(-c_2(N) \rho)$ of the boundary of the space of
admissible polar metrics, where $c_1, c_2$ are strictly positive functions of
$N$. 
\end{theorem}
The constants in the statement of the Theorem above are completely explicit,
and can be sharpened by taking into consideration finer invariants of the
combinatorics of $P$ than the number of vertices. 

In Section \ref{vanish} we give the  proof of the
Vanishing Conjecture for simplices (that is, Milnor's original
conjecture) and then use our description (as given in Section
\ref{boundary}) of the boundary of the set of polar metrics of convex
polytopes to show the Vanishing Conjecture for arbitrary convex polytopes.

The main result of Section \ref{boundary} is
as follows:
\begin{theorem}
Let $P$ lie on the boundary of the space of polar metrics of compact
convex polytopes in $\mathbb{H}^n.$ Then either $P$ has a
\emph{combinatorial} closed geodesic of length $2\pi,$ or $P$ is a
metric suspension.
\end{theorem}
In the above, a combinatorial geodesic is one which is contained in
the $1$-dimensional skeleton of the cell-decomposition of $P$ coming
from a family of polar metrics of degenerating polytopes.

\section{A simple proof for simplices (among other things)}
\label{simplices}

In dimension $2,$ the result follows immediately from Gauss' formula, which
states that area is a linear function of the angles, so we will only discuss
dimensions $3$ or above. 

The simple proof relies on the \emph{Schl\"afli differential equality} (see 
\cite{milnorcol}, which states that in a space of constant curvature $K$ and
dimension $n$ the volumes of a smooth family of polyhedra $P$ satisfy the
differential equation: 

\begin{equation}
\label{schlafli1}
K d V(P) = \dfrac{1}{n-1} \sum_F V_{n-2}(F) d \theta_F,
\end{equation}

where the sum is over all codimension-$2$ faces, $V_{n-2}$ is the
$n-2$ dimensional volume of $F,$ and $\theta_F$ is the dihedral angle
at $F.$ 

Another way of writing the Schl\"afli formula is:

\begin{equation}
\label{schlafli2}
K \dfrac{\partial V(P)}{\partial \theta_F} = V_{n-2}(F).
\end{equation}

This is the form we will use. 

The first observation is that $V_{n-2}(F)$ is bounded by a constant
(dimensional for  $\mathbb{S}^n,$ depending on the number of vertices of
$F$ in $\mathbb{H}^n$ for $n \geq 4.$)  

This immediately shows the continuity of volume for all $\mathbb{S}^n,$ and
for $\mathbb{H}^n,$ whenever $n \geq 4.$ 

We are left with dimension $3.$ All we really need is the result that the partial
derivatives of $V$ with respect to the dihedral angles develop at worst
logarithmic singularities as we approach the frontier of
$\overline{\Omega_{\mathbb{H}^3}}$ -- this result suffices by the following
form of the Sobolev Embedding Theorem (this is \cite[Theorem 7.26]{gilbrtru}): 
\begin{theorem}
\label{sobemb}
Let $\Omega$ be a $C^{0,1}$ domain in $\mathbb{R}^n.$ Then,
\begin{itemize}
\item{(i)} If $k p < n,$ the space $W^{k,p}(\Omega)$ is continuously
imbedded in $L^{p^*}(\Omega),$ where $p^* = np/(n - kp),$ and
compactly imbedded in $L^q(\Omega)$ for any $q < p^*.$ 
\item{(ii)} If $0 \leq m < k - \frac{n}{p} < m+1,$ the space $W^{k,
p}$ is continously embedded in $C^{m, \alpha}(\overline{\Omega}), $
$\alpha=k -n/p - m,$ and compactly embedded in $C^{m,
\beta}(\overline{\Omega})$ for any $\beta < \alpha.$ 
\end{itemize}
\end{theorem}
Here, the Sobolev space $W^{k, p}$ is the space of functions whose
first $k$ (distributional) derivatives are in $L^p.$

In our case, we know that the domain $\Omega$ is bounded, convex
``curvilinear polyhedral'' (hence $C^{0, 1}$) domain,  volume is a
bounded function, and we assume that the gradient grows
logarithmically as we approach the boundary. This implies that $V$ is
in $W^{1, p}$ for \emph{all} $p > 0,$ so we get the following
corollary: 
\begin{corollary}
Volume is in $C^{0, \alpha}(\overline{\Omega})$ for any $\alpha < 1.$
\end{corollary}

The logarithmic growth of diameter of the simplex as a function of the
distance to $\partial \Omega$ can be shown in a completely elementary
way using Eq. \eqref{schlafli2} and elementary reasoning about Gram
matrices, as follows:

Let $G$ be ``angle Gram matrix'' of a simplex $\Delta$, that is,
$G_{ij} = - \cos\theta_{ij},$ where $\theta_{ij}$ is the angle between
the $i$-th and the $j$-th face.  Let $S$ be the matrix whose columns
are the normals to the faces of $\Delta$ (all the computations take
place in Minkowski space, and we use the hyperboloid model of
$\mathbb{H}^n.$ It is immediate that $G = S^t S.$

Let now $W$ be the matrix whose columns are the (possibly scaled)
vertices of $\Delta.$ $W$ satisfies the equation $S^t W = I,$ and to
get the vertices to lie on the hyperboloid $\langle x, x\rangle = -1$
we must rescale in such a way that the squared norms of the columns of
$W$ become $-1.$ Call the scaled matrix $W_s.$ Since the usual
``length'' Gram matrix $G^*$ of $\Delta$S can be written as
$W_s^tW_s,$ and $G^*_{ij} = - \cosh(d(v_i, v_j)),$ a simple
computation using Cramer's rule gives:

\[
\cosh d(v_i, v_j) = \dfrac{c_{ij}}{\sqrt{c_{ii} c_{jj}}},
\]
where $c_{ij}$ is the $ij$-th cofactor of $G.$ (see \cite{muraushi}
for many related results).

It  follows that the distances between the vertices (which are the
lengths of the edges, which are the faces of codimension $2.$) behave
as $|\log\ c_{ii}|.$ Since the cofactors are polynomial in the cosines
of the angles, we are done. 

It should be noted that this argument works \emph{mutatis mutandis}
for \emph{hyperideal} simplices, or simplices with some finite and
some hyperinfinite vertices.. 

\section{Convex polytopes}
\label{convexpoly}
For arbitrary convex polytopes in dimension $n > 3$ (and convex
\emph{spherical} polytopes in all dimensions)  the proof given in
Section \ref{simplices} goes through without change, with the one
proviso that it is not currently known whether the volume of a
polytope is determined up to congruence by its dihedral angles. Such a
uniqueness result \emph{is} conjectured (indeed, it is conjectured
that a polytope is determined up to congruence by the dihedral
angles), and is easy to prove for \emph{simple} polytopes -- those
with simplicial links of vertices -- this follows in arbitrary
dimension from the corresponding result in $3$ dimensions
(\cite{thes,thes0}).  The uniqueness 
issue can be finessed (in dimension $3,$ at least) by using the results
of \cite{thes,thes0}:
\begin{theorem}[\cite{thes,thes0}]
\label{annelk}
A metric space $(M, g)$ homeomorphic to $\mathbb{S}^2$ can arise as
the Gaussian image $G(P)$ of a compact convex polyhedron $P$ in
$\mathbb{H}^3$ if and only if the following conditions hold: 
\begin{itemize}
\item{(a)} The metric $g$ has constant curvature $1$ away from a
finite collection of cone points. 
\item{(b)} The cone angle at each $c_i$ is greater than $2\pi.$
\item{(c)} The lengths of closed geodesics of $(M, g)$ are all
strictly greater than $2\pi.$ 
\end{itemize}
\end{theorem}
The space of
admissible metrics $\Omega_P$ (as per Theorem \ref{annelk}) is
parametrized by the exterior dihedral angles (the cell decomposition
dual to that of $P$ gives a triangulation of the Gaussian image, and
the (exterior) dihedral angles are the lengths of edges of the
triangulation.) Theorems \ref{degenthm},\ref{quasig} immediately imply
the following: 

\begin{theorem}
\label{loggrowth}
There exists a constant $L_0,$ such that 
the maximal length $\ell_P$ of an edge of $P$ is bounded as follows:
\[
\ell_P \leq \max(L_0,  - 2N \log(d(P, \partial \overline{\Omega}_P)/12N)),
\]
where $N$ is the number of vertices of $P.$
\end{theorem}
\begin{proof}
Assume the contrary. Then, there exists a sequence of polyhedra $P_1,
\dots, P_n, \dots$ with diameter $\rho(P_i) \geq \ell){P_i}$ going to
infnity, which are farther than $12 N  \exp(-\rho/2N).$ By choosing a
subsequence, we may assume that there is a \emph{fixed} cycle of faces
$F_1, \dots, F_k$ of $P,$ such that the sum of dihedral angles along
the edges $e_i = F_i \cap F_{i+1}$ is smaller than $2\pi + 12 N
\exp(-\rho/2N),$ (by Theorem \ref{degenthm}) and which are a
$4N\exp(-2\rho)$ quasigeodesic (by Theorem \ref{quasig}). Since the
limit point of the $P_i$ is not in $\Omega_P$ (by Theorem
\ref{annelk}), 
the result follow.
\end{proof}
The following corollary is immediate (by Schl\"afli, see Section
\ref{simplices}): 
\begin{corollary}
\label{lpest}
The volume is in $W^{1, p}(\Omega_p)$ for \emph{all} $p > 0.$
\end{corollary}
We now have almost enough to show that volume extends to
$\overline{\Omega}_p,$ except for the slight matter of not having the
required (by Theorem \ref{sobemb}) regularity result for
$\partial{\Omega}_P.$ Such a result seems quite non-trivial, since the
length of the shortest closed geodesic is a rather badly behaved
quantity, but the results of Section \ref{boundary} show that things
are well enough behaved. 

\section{Degeneration estimates}
\label{degennes}
The results of this section are a quantitative version of the results
of the compactness results of \cite{thes,thes0}.  
First, some key lemmas. The general setup will be as follows: $L$ is a
geodesic in $\mathbb{H}^3,$  $t$ is a real number (generally large)
and $P, P^-, P^+$ are three planes, all orthogonal to $L,$ and such
that $d(P, P^-) = d(P, P^+) = t,$ and $d(P^-, P^+) = 2 t.$ We denote
$x_0 = L \cap P.$ 

In the sequel, we use the hyperboloid model of $\mathbb{H}^3,$ where
$\mathbb{H}^3$ is represented by the set $\langle x, x \rangle = -1;$
$x_0 > 0,$ in the $\mathbb{R}^4$ equipped with the scalar product
$\langle x, y \rangle = - x_1 y_1 + \sum_{i=2}^4 x_i y_i.$ The reader
is referred to \cite{thurston97} (as well as \cite{thes}) for the details
(which will be used below). 

Returning back to our setup, we can assume, without loss of generality, that
\[
x_0=\begin{pmatrix}1\\0\\0\\0\end{pmatrix},
\]
that 
\[
P^\perp = \begin{pmatrix}0\\1\\0\\0\end{pmatrix},
\]
and hence, that $P^+ = \phi(t) P,$ while $P^- = \phi(-t) P,$ where
\[
\phi(r) = \begin{pmatrix} \cosh(r) & \sinh(r) & 0 & 0\\
                                           \sinh(r) & \cosh(r) & 0 & 0\\
                                           0 & 0 &1 &0\\
                                           0 & 0 & 0 &1
                                           \end{pmatrix}.
                                           \]
                                           
Since $\phi(r)$ is symmetric, it follows that \[P^{+\perp} = \phi(t)
P^\perp = \begin{pmatrix} \cosh(t) \\ \sinh(t) \\ 0
\\0\end{pmatrix},\] while  
\[P^{-\perp} = \phi(t) P^\perp = \begin{pmatrix} \cosh(t) \\ -\sinh(t)
\\ 0 \\0\end{pmatrix},\] 
\begin{lemma}
\label{anglemma}
 Let $Q$ be a plane in $\mathbb{H}^3$ which intersects both $P^-$ and
$P^+.$ Then, there exists $t_0,$ such that $Q$ intersects $P,$ and the
cosine of the  angle $\alpha$ of intersection satisfies
$|\cos(\alpha)| <3 e^{-t},$ as long as $t > t_0.$ The number $t_0$ can
be picked \emph{independently} of $Q.$ 
\end{lemma}
 
\begin{proof}
Let the unit normal $Q^\perp$ to $Q$ be
$Q^\perp=\begin{pmatrix}a\\b\\c\\d\end{pmatrix}.$   
Since two planes intersect if and only if the scalar product of their
unit normals is less than $1$ in absolute value, we have, from the
hypotheses of the lemma and the description of the unit normals to
$P^-$ and $P^+$ above that: 
\begin{gather}
|a \cosh(t) + b \sinh(t)| < 1 \\
|a \cosh(t) - b \sinh(t)| < 1.
\end{gather}
Squaring the two inequalities, and adding them together we obtain:
\[
a^2 \cosh^2(t) + b^2 \sinh^2(t) < 1.
\]
Since, under the hypotheses of the lemma, $\min(\cosh(t), \sinh(t)) >
e^t/3,$ it follows that 
\[
a^2 + b^2 < 3/e^t,
\]
and so $\max(a, b) < 3e^{-t}.$
Now, the cosine of the angle between $Q$ and $P$ equals $\langle
Q^\perp, P^\perp\rangle = b,$ so the result follows. 
\end{proof}

\begin{remark}
The constant $3$ is far from sharp (especially for larger $t$).
\end{remark}

\begin{lemma} 
\label{distlem}
There exists a $t_0,$ such that if  $M$ is a line in $\mathbb{H}^3$
which intersects both $P^-$ and $P^+,$ then $M$ intersects $P,$ and
$\cosh(d(P\cap M, x_0)) < 4 e^{-2t}+1, $ as long as $t > t_0.$ 
\end{lemma}

\begin{proof}
Assume that 
$M\cap P^+ = \phi(t)  p_1,$ and $M\cap P^- = \phi(-t) p_1,$ where
$p_{1, 2} \in P.$ (This is always possible, since $P^+ = \phi(t) P,$
$P^- = \phi(-t).$) The intersection of $M$ with $P$ is then given by  
\[
M\cap P = \dfrac{x (M \cap P^+) + y(M \cap P^-)}{\|x (M \cap P^+) +
y(M \cap P^-)\|}, 
\]
where $x$ and $y$ are chosen so that the linear combination is
actually in $P,$ or, in other words, the second coordinate of the
linear combination vanishes. We abuse notation above by writing 
$\|Z\| = \sqrt{-\langle Z, Z\rangle}.$

Let us now compute. Set (for $i=1, 2$)
\[
p_i = \begin{pmatrix} a_i\\0\\c_i\\d_i\end{pmatrix}.
\]
It follows that 
\[
M\cap P^+ = \begin{pmatrix}a_1 \cosh(t) \\ a_1 \sinh(t)\\ c_1 \\
d_1\end{pmatrix}. 
\]
while
\[
M \cap P^- = \begin{pmatrix}a_2 \cosh(t) \\ - a_2 \sinh(t) \\ c_2
\\d_2 \end{pmatrix}. 
\]
It follows that we can choose $x = 1/(2a_1),$ $y = 1/(2a_2),$
so that
\[
m = x M\cap P^+ + y M\cap P^- = \begin{pmatrix} \cosh(t) \\0 \\
\frac{1}{2}( c_1/a_1 + c_2/a_2) \\ \frac{1}{2} ( d_1/a_1 + d_2/a_2)
\end{pmatrix}. 
\]
It follows that 

\begin{multline}
- \cosh(d(M \cap P, x_0)) = \left\langle \dfrac{m}{\|m\|},
x_0\right\rangle = \\ -\dfrac{\cosh(t)}{\sqrt{\cosh^2(t) -
1/4\left((c_1/a_1 + c_2/a_2)^2 + (d_1/a_1 + d_2/a_2)^2\right)}}. 
\end{multline}

Since $c_i^2 + d_i^2 + 1 = a_i^2,$ for $i=1,2$ it follows that
$|c_i/a_i| < 1,$ and similarly $|d_i/a_i| < 1,$ so that  
\[
\cosh^2(t) \geq \cosh^2(t) - 
 1/4\left((c_1/a_1 + c_2/a_2)^2 + (d_1/a_1 + d_2/a_2)^2\right) > \cosh^2(t) -2.
 \]
 It follows that 
 \[
\cosh(d(M \cap P, x_0))\leq \dfrac{1}{\sqrt{1-2/\cosh^2(t)}},
 \]
 and the assertion of the lemma follows by elementary calculus.
\end{proof}

\begin{lemma}
\label{spherical}
Let $T$ be a spherical triangle with sides $A,B,C$ and (opposite)
angles $\alpha, \beta, \gamma.$ Suppose that $|\cos(\beta)| < \epsilon
\ll 1,$ $| \cos(\gamma)| < \epsilon \ll 1.$ Then  
$|\alpha - A| < 2 \epsilon|.$
\end{lemma}

\begin{proof} The spherical Law of Cosines states that:
\[
\cos(A) = \dfrac{\cos(\alpha)+ \cos(\beta) \cos(\gamma)}{\sin(\beta)
\sin(\gamma)}. 
\]
It follows that 
\[
\cos(A) - 2\epsilon^2 \leq \cos(A)(1-\epsilon^2) - \epsilon^2  \leq
\cos(\alpha) \leq \cos(A) + \epsilon^2. 
\]
The assertion of the lemma follows immediately.
\end{proof}

\begin{corollary}
\label{spherecor}
Let $F_1$ and $F_2$ be two planes intersecting at a dihedral angle
$\alpha,$ with both $F_1$ and $F_2$ intersecting a third plane $P,$ at
angles whose cosines are smaller than $\epsilon.$ Let $A$ be the angle
between $F_1 \cap P$ and $F_2 \cap P.$ Then $|\alpha - A|< 2\epsilon.$

\end{corollary}

\begin{proof}
Apply Lemma \ref{spherical} to the link of the point $F_1 \cap F_2 \cap P.$
\end{proof}

\begin{lemma}
\label{circleest}
Let $V$ be a convex polygon in the hyperbolic plane $\mathbb{H}^2,$
such that all the vertices of $V$ lie within a distance $r$ of a
certain point $O.$ Then, the sum of the exterior angles of $V$ is
smaller than $2\pi \cosh(r).$ 
\end{lemma}

\begin{proof}
The area of a disk of radius $r$ in $\mathbb{H}^2$ equals $4\pi
\sinh^2(r/2) = 2\pi (\cosh(r) -1)$ (see \cite{vinberg}).  
Since $V$ is contained in such a disk, its area is at most
$2\pi(\cosh(r) -1),$ and since the area of $V$ equals the difference
between the sum of the exterior angles and $2\pi,$ the statement of
the lemma follows. 
\end{proof}

Now we are ready to show the following:
\begin{theorem}
\label{degenthm}
Let $X$ be a convex polyhedron with $N$ vertices  in $\mathbb{H}^3$ of
diameter $\rho \gg 1.$ Then, there exists a cyclic sequence of faces
$F_1, \dotsc, F_k = F_1,$ with $F_i$ sharing an edge $e_i$  with
$F_{i+1}$ (indices taken $\mod k$) so that the sum of exterior
dihedral angles at $e_1, \dotsc, e_k$ is smaller than $2\pi + 12 N
\exp(-\rho/2N).$ 
\end{theorem}

\begin{proof}
Take a diameter $D$ of $X$  of length $\rho,$ place points $p_1,
\dotsc, p_N$ equally spaced on $D.$ By the pigeonhole principle, one
of the segments $p_ip_{i+1}$ contains no vertices of $X.$ Let $x_0$ be
the midpoint of the segment $p_ip_{i+1}.$ Construct planes orthogonal
to $D$ at $x_0$ ($P$) and $p_i$ ($P^-$), and at $p_{i+1}$ ($P^+$). Let
$t=\rho/(2N).$ The portion of $X$ contained between $P^-$ and $P^+$ is
a polyhedral cylinder, consisting of faces $F_1, \dots, F_k.$ By Lemma
\ref{distlem}, the intersection of $X$ with $P$ is a polygon
$\mathcal{P},$ whose sum of exterior angles is at most
$2\pi(4\exp(-2t) + 1),$ and so by Corollary \ref{spherecor}, combined
with Lemma \ref{anglemma}, the sum of the dihedral angles
corresponding to pairs $F_i F_{i+1}$ is at most $2\pi(4 \exp(-2t) + 1)
+ 6 k \exp(-t).$ Since $k$ is no greater than the number of faces of
$X,$ which, in turn, is at most $2N -4.$ 
\end{proof}

\begin{theorem}
\label{quasig}
With notation as in Theorem \ref{degenthm}, the faces $F_1, \dotsc,
F_k$ form a  curve in the Gaussian image of $X$ with geodesic
curvature not exceeding $3k\exp(-\rho/N)).$ 
\end{theorem}
\begin{remark} The reader is referred to \cite{thes0,thes} for a more
thorough discussion of geodesics on spherical cone manifold, but
suffice it to say that the contribution of the face $F_i$ to the
geodesic curvature is $0$ if the two edges are (hyper)parallel, and
equal to the angle of intersection if they intersect. 
\end{remark}
\begin{proof}
Let $e_1$ and $e_2$ be the two edges of $F.$ If $e_1$ and $e_2$ do not
intersect, there is nothing to prove (by the remark above. If they do
intersect at a point $C,$ note that $C$ is at a distance at least
$\rho/2N$ from $x_0,$ while the intersections $A$ and $B$ of $e_1$ and
$e_2$ with $P$ are at most  
$\acosh(4\exp(\-rho/N)+1) \approx \sqrt{8} \exp(-\rho/2N)$ away from
$x_0,$ and so at most (for large $\rho$) $6\exp(-\rho/2N)$ away from
each other. We will only use the (much cruder) estimate $\cosh(AB)
\leq 2.$ Now, apply the hyperbolic law of cosines to the triangle
$ABC,$ to get: 

\begin{multline}
1\geq \cos(\gamma) = \dfrac{-\cosh(AB) + \cosh(AC)\cosh(BC)}{\sinh(AC)
\sinh(BC)} \geq \\ 
1- \dfrac{\cosh(AB)}{\sinh(AC)\sinh(BC)} \geq 1-
\dfrac{2}{\sinh(AC)\sinh(BC)} \geq 1 - 8 \exp(-\rho/2N). 
\end{multline}
The estimate now follows.

\end{proof} 

\begin{remark}
The argument above is easily modified to show that the curve dual to
$F_1, \dotsc, F_k$ has small geodesic curvature viewed as a curve in
$\mathbb{S}_1^3,$ and not just in $X^*.$ 
\end{remark}

\section{The boundary of the space of polar metrics of convex
polytopes.} 
\label{boundary}
Consider a sequence of degenerating polytopes. We have two
possibilities: the diameter stays bounded or it does not. If the
diameter does not stay bounded, then the results of Section
\ref{degennes} indicate that one can pick a subsequence in such a way
that the length of a (quasi)-geodesic in the dual $1$-skeleton
converges to $2\pi,$ while the quasi-geodesic itself converges to a
dual $1$-skeleton geodesic. The other possibility is that the
polytopes degenerate while the diameter is bounded. In this case there
are the following possibilities:

First, the diameter goes to 0. In this case, it is clear that the
polar is a round sphere. 

Secondly, the diameter stays bounded away from 0, but the limit is
$1$-dimensional. In this case the polar metric is still a round
sphere.

Thirdly, the limit may be $2$-dimensional (a doubled polygon). In this
case the polar is a 
metric suspension with two cone points with curvature equal to the
area of the (doubled) polygon.

In higher dimensions the analysis is the same, though the number of
suspension possibilities increases.

\section{The Vanishing Conjecture}
\label{vanish}

We will first need the following observation:
\begin{lemma}
\label{ballcomp}
The set of (hyper)planes intersecting a fixed ball in $\mathbb{H}^n$
is compact. 
\end{lemma}
\begin{proof}
There are a number of arguments, the simplest of which would appear to
be that the set of planes going through a fixed point in
$\mathbb{H}^n$ is compact (being in one-to-one correspondence with the
unit sphere $\mathbb{S}^{n-1}$) and then identifying the set of planes
intersecting a ball $B$ with a quotient of $\mathbb{S}^{n-1} \times B.$
\end{proof}

We will actually need the following:
\begin{corollary}
A sequence of polytopes all faces of which intersect a fixed ball
contains a convergent subsequence.
\end{corollary}
\begin{proof}
Immediate by compactness.
\end{proof}

In order to deal with the vanishing conjecture for simplices, we now
make the following:
\begin{observation}
\label{incircle}
There exists a universal constant $K$ such that for any triangle
$T\subset \mathbb{H}^2,$ there exists a disk of radius $K$
intersecting all of the sides of $T.$
\end{observation}
The observation can be rephrased as saying that the hyperbolic plane
is Gromov-hyperbolic. The constant $K$ can be chosen to be $\log 2/2.$
\begin{proof}
Since every triangle is contained in an ideal triangle, it is enough
to show the result for the ideal triangle. There, the result follows
by construction.
\end{proof}

\begin{corollary}
There exists a universal constant $K$ such that for any simplex
$T\subset \mathbb{H}^n,$ there exists a ball of radius $K$
intersecting all of the faces of $T.$ 

\end{corollary}

\begin{proof}
By induction on dimension. Pick any face $F$ of the simplex $T\in
\mathbb{H}^n.$ By induction, there is an $n-1$-dimensional ball
of radius $K$ which intersects all of the faces of $F$, and thus all of
the faces of $T.$ 
\end{proof}

Observation \ref{incircle} shows that any sequence of simplices
contains a convergent subsequence, and hence the volume of a sequence
of simplices with degenerating dihedral angles is the volume of an
actual simplex $T_\infty$ in $\mathbb{H}^n.$ The only way that volume
could be equal to $0$ is if $T_\infty$ is degenerate (that is, lower
dimensional). It is easy to see that the dihedral angles of $T_\infty$
then lie in the closure of the set of angles of Euclidean simplices.

To show the Vanishing Conjecture for an arbitrary sequence of
polytopes, we consider two possibilities. The first is that that all
the faces of the polytopes of the (sub)sequence intersect a fixed
ball. This case is the same as the case of the simplex consider above,
and there is nothing left to prove. 

For the other possibility, we will first need the following:
\begin{lemma}
\label{seplemma}
Let $T$ be a simplex in $\mathbb{H}^n,$ and $B$ a ball intersecting all the faces of
$T.$ Let $P$ be a plane which does not intersect $B.$ Then at least
$2$ vertices of $T$ lie on the same side of $P$ as $B.$
\end{lemma}
\begin{proof}
Suppose not. Then at least $n$ vertices of $T$ are separated from $B$
by $P,$ and hence so is their convex hull, which is then a face of $T$
not intersecting $B,$ contradicting the hypothesis.
\end{proof}

\begin{corollary}
\label{sepcol}
Let $T_1$ and $T_2$ be two simplices in $\mathbb{H}^n,$ let $B_1$ and
$B_2$ be balls intresecting all the faces  of $T_1$ and $T_2,$ respectively, and let $P$
be a hyperplane such that $B_1$ and $B_2$ are on different sides of
$P,$ and which \emph{does not} contain $T_1 \cap T_2$. Then there are
at least $2$ vertices of $T_1$ on one side of $P$ and at least $2$ vertices of $T_2$
on the other side.
\end{corollary}
\begin{proof}
Follows immediately from Lemma \ref{seplemma}
\end{proof}

Let us now assume that there is no ball which all the faces
intersect. Let us assume, for convenience, that all the faces of the
polytopes in the sequence are simplicial (if not, we can always
triangulate them, with the additional dihedral angles equal to $\pi.$
For each face $F_i$ we have the ball $B_i$ which interesects all of
its faces and there must
be a pair of adjacent faces $F_i, F_j$ such that the $B_i$ and $B_j$ are far
apart. Let $E_{ij}$ be $F_i \cap F_j,$ and there must be a cycle of
faces $f_1 = F_i, f_2 = F_j, f_3, \dotsc, f_n = f_1$ which give a dual
quasi-geodesic of length close to $2\pi$ and a corresponding plane $P$
(as in Section \ref{degennes}), By the lemma,
the set of vertices of our polytope is separated by $P$ into two sets,
the cardinality of each of which is at least $2,$ 
 and the limiting
object is the disjoint union of two limits, one on each side of $P,$
and the limiting volume is the sum of the two volumes. One can then
induct on the number of vertices to show that both halves are
degenerate, and hence so is the limit.

\bibliographystyle{plain}
\bibliography{curves,rivin,opt}

\begin{thebibliography}{10}

\bibitem{gilbrtru}
David Gilbarg and Neil Trudinger.
\newblock {\em Elliptic Partial Differential Equations of Second Order}.
\newblock Classics in Mathematics. Springer Verlag, Berlin, New York, 1998.

\bibitem{kneser36}
Hellmuth Kneser.
\newblock Der {S}implexinhalt in der nichteuclidischen {G}eometrie.
\newblock {\em Deutsche Mathematik}, 1:337--340, 1936.

\bibitem{luovol}
Feng Luo.
\newblock Continuity of volume of simplices in classical geometry.
\newblock Technical Report math.GT/0412208, arxiv.org, 2004.

\bibitem{milnorcol}
John~W. Milnor.
\newblock {\em The {S}chl{\"a}fli differential equality}, volume~1.
\newblock Publish or Perish, Houston, Texas, 1994.

\bibitem{muraushi}
Jun Murakami and Akira Ushijima.
\newblock A volume formula for hyperbolic simplices in terms of edge lengths.
\newblock Technical Report math.MG/0402087, arxiv.org, 2004.

\bibitem{thes0}
Igor Rivin.
\newblock {\em On the Geometry of Convex Polyhedra in Hyperbolic 3-Space}.
\newblock PhD thesis, Princeton University, July 1986.

\bibitem{rivvol}
Igor Rivin.
\newblock Continuity of volumes -- on a generalization of a conjecture of
  {J.W.M}ilnor.
\newblock Technical Report math.GT/0502543, arxiv.org, 2005.

\bibitem{thes}
Igor Rivin and C.D.Hodgson.
\newblock A characterization of compact convex polyhedra in hyperbolic 3-space.
\newblock {\em Inventiones Mathematicae}, pages 77--111, January 1993.
\newblock Corrigendum, vol 117, page 359.

\bibitem{thurston97}
William~P. Thurston.
\newblock {\em Three-dimensional geometry and topology, vol. 1}.
\newblock Number~35 in Princeton Mathematical Series. Princeton University
  Press, Princeton, New Jersey, 1997.

\bibitem{vinberg}
E.~B. Vinberg.
\newblock {\em Geometry {II}}, volume~29 of {\em Encyclopaedia of
  {M}athematical {S}ciences}.
\newblock Springer Verlag, Berlin-Heidelberg-New York, 1993.

\end{thebibliography}
\end{document}